\documentclass[12 pt]{article}
\usepackage{verbatim}
\usepackage{amsmath}
\usepackage{latexsym}
\usepackage{amssymb}
\usepackage{amscd}
\usepackage{amsthm}
\usepackage{amssymb}  
\usepackage{latexsym} 
\usepackage{comment}
\usepackage{graphicx}
\usepackage{graphics}
\usepackage{epsfig}
\newtheorem{Thm}{Theorem}[section]

\newtheorem{Lem}[Thm]{Lemma}

\newcommand{\Z}{{\mathbb Z}}
\newcommand{\Q}{{\mathbb Q}}
\newcommand{\R}{{\mathbb R}}

\begin{document}

\openup5pt

\title{Mordell-Weil Problem for Cubic Surfaces, Numerical Evidence}

\author
{Bogdan G. Vioreanu
\thanks{Yale University, New Haven, CT, USA, 
E-mail: bogdan.vioreanu@yale.edu}}
\date{}
\maketitle

\begin{abstract} 
Let $V$ be a plane smooth cubic curve over a finitely generated field $k$. The Mordell-Weil theorem
for $V$ states that there is a finite subset $P \subset V(k)$ such that the whole $V(k)$ can be obtained 
from $P$ by drawing secants and tangents through pairs of previously constructed points and consecutively
adding their new intersection points with $V.$ 
In this paper we present numerical data regarding the analogous statement for cubic surfaces.
For the surfaces examined, we also test Manin's conjecture relating the 
asymptotics of rational points of bounded height on a Fano variety with the rank of the
Picard group of the surface.  
\end{abstract}

\section{Introduction}

Let $V$ be a smooth cubic surface over a field $k$ in $\mathbb{P}^3.$
If $x,y,z \in V(k)$ are three points (with multiplicities) lying on a line in $\mathbb{P}^3$ not belonging
to $V,$ we write $ x = y \circ z $. Thus $\circ$ is a partial and multivalued composition law on $V(k)$. 
Note that $x \circ x$ is defined as the set of points in the intersection of $V(k)$ with the tangent plane at $x$.
If $x$ does not lie on a line, this is a cubic curve $C(x)$ with double point $x \in V(k)$. 
This whole set must be considered as the domain
of the multivalued expression $x \circ x$, because geometrically all its points can be obtained by drawing 
tangents with $k$-rational direction to $x$. This means that
an important source for generating new rational points on the
cubic surface will be doubling the points that were already generated.
The analogue of the Mordell-Weil theorem for cubic surfaces states that 
$(V(k), \circ)$ is finitely generated, i.e., there is 
a finite subset $P \subset V(k)$ such that the whole $V(k)$ can be obtained from $P$ by drawing
secants and tangent planes through pairs of (not necessarily distinct) previously constructed points,
and consecutively adding their new intersection points with $V.$ By drawing secants we can add only one 
rational point to $P,$ while tangent sections give us an infinite number of points that can be generated, 
by the note above. 
For a more thorough discussion of various versions of finite generation cf. \cite{K&M}.
Note that, by Theorem 11.7 of \cite{Man1},  
finite generation of $(V(k), \circ)$ implies that the universal quasi-group of
$(V(k), \circ)$, as defined in \cite{Man1}, chapter II,  is finite 
and has $2^n 3^m$ elements for some $n,m \in \Z_{\geq0}$.    

In the following, we present the procedure we used to test 
whether $(V(\Q), \circ)$ is finitely generated, and the results we obtained for thirteen diagonal
cubic surfaces, six of them having the rank of their Picard group equal to 1, and seven of them mentioned in \cite{P&T}, illustrating the cases of surfaces with ranks $2$ and $3$ of the Picard group. We also bring 
numerical evidence supporting Manin's conjecture for the asymptotics of rational points
of bounded height on a Fano variety. Note that John Slater and Sir Peter Swinnerton-Dyer
have proved in \cite{S&S-D} a one-sided estimate for the conjecture in the case when $V$ 
contains two rational skew lines. 
All the computations were done using the Magma computer algebra system (cf. \cite{Magma}.)

\section{Description of the procedure}

Let $a x^3 + b y^3 + c z^3 + d u^3 = 0$, where $a, b, c, d$ are nonzero integers, be a
diagonal cubic surface. Using a program due to Dan Bernstein (see \cite{Ber}), we find all rational 
points on this surface up to height $H = 10^5$ or $H = 1.5 \cdot 10^5$, 
where the height of a rational point $P = (x:y:z:u)$,
with $x,y,z,u \in \Z$ and $\gcd(x,y,z,u) = 1$ is defined as 
\[ h_{max}(P) := \max \{\left|x\right|, \left|y\right|, \left|z\right|, \left|u\right| \} \,. \]
We consider also another height function $h_{sum}: V(\Q) \longrightarrow \R_+$ defined by 
\[ h_{sum}(P) := \left|x\right| + \left|y\right| + \left|z\right| + \left|u\right| \,. \]
Note that a rational point $P$ can be uniquely written in the above form up to a sign change 
of the coordinates. So, if we assume, in addition, that the first nonzero coordinate of $P$
is positive, then there is a unique such 'canonical' form corresponding to each point $P.$  
We order the rational points by increasing $h_{sum}$. If there are two or more points having
the same height $h_{sum}$, then we order them lexicographically according to their coordinates
in the canonical form. This defines a total order on the set of rational points. 
We will write $P < Q$ if $P$ precedes $Q$ in the 
sorted list, and use the number of a point in this list as its name. We will also 
refer to this number as the {\em index} of a rational point. 

We will use the $h_{max}$ height function only to study the 
asymptotics of the number of rational points on a cubic surface, 
while for the ordering of the points and in the implementation
of the main function we will use $h_{sum}$.     
  
For testing whether a given set of rational points is generating, we use the procedure
{\em Test Generating Set (TGS)}, which is described below. 

The procedure implements essentially a descent method. 
Given an index bound $n$ and a set of points $GeneratedSet$ that is presumably generating, we perform the following
iterative process. In one iteration of  loop, we consider all points in the range $\{1, \ldots, n\}$ that are not in $GeneratedSet$ and test whether they can be decomposed as $x \circ y$, with $x,y \in GeneratedSet$. 
Every point that can be decomposed in such a way is added to the $GeneratedSet$ and at the end of the loop, 
the procedure is reiterated. As now $GeneratedSet$ is bigger, there may be additional points
in the range $\{1, \ldots, n\}$ that can be generated because we can choose the 
points $x,y$ for a possible decomposition from a bigger
set. The procedure is repeated until $GeneratedSet$ stabilizes, 
i.e., until some iteration of the loop does not add any new 
points to the $GeneratedSet$.

In order to avoid repeating some operations of composing points, 
we use the additional variables $OldGeneratedSet$, $JustAdded$ and $Decomp$.  
$OldGeneratedSet$ stores the value of $GeneratedSet$ at the beginning of the iteration of the loop. At the end
of the preceding loop, a number of points will have been added to $GeneratedSet$. These points are stored in the 
set variable $JustAdded$. During an iteration of the loop, we store in $Decomp$ decompositions 
of the type $i = j \circ k$, with $i,j,k \leq n, \, i,k \notin GeneratedSet$ and $j \in GeneratedSet$.
These are the only decompositions that we could further use. Indeed, if, at some point, $k$ was added to
$GeneratedSet$, then by searching in $Decomp$, we would find the decomposition $j \circ k$ of $i$ and we would
add $i$ to $GeneratedSet$ without performing any composition of points (which requires multiplications, 
so is computationally expensive) because we know, by the way we constructed $Decomp$, that 
$j \in GeneratedSet$ already. 

Receiving as input the parameters $GeneratedSet$ (a set of points in $V(\Q)$ that is 
assumed to be generating), and $n$ (the index bound for the points used in the 
decompositions), the $TGS$ procedure does the following:       
 
\begin{enumerate}\addtolength{\itemsep}{0mm}
	\item Set $Decomp = \emptyset$, $OldGeneratedSet= \emptyset$. 
	\item Set $JustAdded = GeneratedSet \setminus OldGeneratedSet$,\\	
				\indent \hspace{6mm} $OldGeneratedSet = GeneratedSet$.
		\label{loop} 
	\item If 	$JustAdded = \emptyset$, return $GeneratedSet$.		
	\item For every point $i \in \{1,2, \ldots, n\} \setminus GeneratedSet$ do:\\     
        \indent \hspace{8mm} search in $Decomp$ for decompositions of $i$ 
        											as $x \circ y$ with $y \in JustAdded$\\
        \indent \hspace{8mm} if such a decomposition exists, add $i$ to $GeneratedSet$\\ 
				\indent \hspace{8mm} else for every point $j$ in $JustAdded$ do:\\
					\indent \hspace{16mm} $k = i \circ j$\\
					\indent \hspace{16mm} if $k \in JustAdded$\\
						\indent \hspace{24mm} add $i$ to $GeneratedSet$\\
						\indent \hspace{24mm} break\\
					\indent \hspace{16mm} else if $k \leq n$ add the decomposition $(j \circ k)$ of $i$ to $Decomp$\\
					\indent \hspace{16mm} end for\\
				\indent \hspace{8mm} end for.
	\item Go to step \ref{loop}. 
\end{enumerate}

Let us explain in more detail the way the algorithm works. 
Suppose that an iteration of the outer loop has just finished, and we are in step \ref{loop}.   
We set $JustAdded = GeneratedSet \setminus OldGeneratedSet$ and test whether this is the empty set. If this is so,
then during the last iteration we could not generate any new points, so the maximum set of points
that can be generated is the current $GeneratedSet$. If $JustAdded$ is not empty, then during the last iteration we found a number of new points that could be generated and added them to $GeneratedSet$
(these are the elements of $JustAdded$), so there is hope
of generating other points. We consider a point $i \notin GeneratedSet$. Since we have already tested during the 
previous iteration whether we could decompose $i$ as $x \circ y$, with $x, y \in OldGeneratedSet$, all we have to 
check now is whether we can write $i = x \circ y$ for $x \in JustAdded$ and either $y \in OldGeneratedSet$ or
$y \in JustAdded$. At the previous iterations of the loop all compositions of $i$ with points in $OldGeneratedSet$ 
that could further be used (i.e., compositions whose result is not bigger than $n$) were stored in $Decomp$, so       
we can check for the first possibility by searching in the vector $Decomp$.
Since by construction we only store in $Decomp$ decompositions of the type $x \circ y$, 
with $x \in GeneratedSet$, all we have to check in the beginning of step $4$ is whether
$y \in JustAdded$ - we are sure that $x \in GeneratedSet$. 
In order to check for the second possibility, we have to compose $i$ with every point $j \in JustAdded$. If the 
result $k$ of the composition is in $JustAdded$, then we can write $x$ as a composition of two points in
$JustAdded$, so we add $i$ to $GeneratedSet$. If the result $k \notin JustAdded$, but    
could be further used (i.e., $k \leq n$), then we store the corresponding
decomposition $j \circ k$ of $i$ in $Decomp$. The 'out of bounds' compositions, i.e., such that $i \circ j > n$,
are implicitly remembered in the process (in the sense that they are done only once.)       

Using the vector $Decomp$ of course implies a tradeoff between space and speed, but we considered 
the latter to be more important. Even with $Decomp$, the computations for $TGS$ for bounds $n$ in the range
of $10^5$ last for several days and sometimes even weeks on an Intel Pentium IV processor with 2.26 GHz.     

Before we proceed with the presentation of the results, let us as provide an estimate of the 
height of the composition of two rational points. Here by $h$ we mean either $h_{max}$ or
$h_{sum}$ since the estimation of the asymptotics does not depend on the choice of the height
function.   

\begin {Lem}
Let $V: a x^3 + b y^3 + c z^3 + d u^3 = 0$ be a diagonal cubic surface, where $a,b,c,d$ are nonzero
integers, and let $K := max \{ \left| a \right|, \left| b \right|, \left| c \right|, \left| d \right| \} $.
If $A_1$ and $A_2$ are two distinct points in $V(\Q)$ that do not lie on a line in $V$, then
\[ h(A_1 \circ A_2) = O(K \cdot max \{ h(A_1), h(A_2) \}^2 \cdot min \{ h(A_1), h(A_2) \}^2 ) \,. \]  
\end{Lem}

$\textit {Proof:}$ Let $A_1 = (x_1 : y_1 : z_1 : u_1)$, $A_2 = (x_2 : y_2 : z_2 : u_2)$ be in
canonical form. Then one can check that 
\[ A_1 \circ A_2 = (\alpha x_1 - \beta x_2 : \alpha y_1 - \beta y_2 : 
										\alpha z_1 - \beta z_2 : \alpha u_1 - \beta u_2)  \,, \]
where 
\[ \alpha = a x_1 x_2^2 + b y_1 y_2^2 + c z_1 z_2^2 + d u_1 u_2^2 \in \Z  \,,\]
\[ \beta  = a x_1^2 x_2 + b y_1^2 y_2 + c z_1^2 z_2 + d u_1^2 u_2 \in \Z  \,.\]
Since the above coordinates of $A_1 \circ A_2$ are integers, the conclusion follows.
This upper bound cannot be improved because, in most cases, the formula given represents
$A_1 \circ A_2$ in canonical form (up to a sign change of the coordinates).

Concerning the doubling of points, if $A \in V(Q)$ is a rational point not lying on a line in $V$, then there 
is no upper bound for the height of the points in $A \circ A$ (since there are infinitely
many such points). On the other hand, there can be many points of small height in 
$A \circ A$, especially if $A$ has small height.    

\section{Results}
Listed below are the thirteen diagonal cubic surfaces that were tested for finite generation, ordered
according to the ranks of their Picard groups:\\

Rank 1 of the Picard group:
\begin{enumerate}\addtolength{\itemsep}{1mm}
\item $x^3 + 2y^3 + 3z^3 + 4u^3 = 0$.		
\label{RS}
\item $x^3 + 2y^3 + 3z^3 + 5u^3 = 0$.		
\label{Second_RS}
\item $17x^3 + 18y^3 + 19z^3 + 20u^3 = 0$.		
\label{Big_RS}
\item $4x^3 + 5y^3 + 6z^3 + 7u^3 = 0$.		
\label{New_RS}
\item $9x^3 + 10y^3 + 11z^3 + 12u^3 = 0$.		
\label{New_Big_RS}
\item $x^3 + 5y^3 + 6z^3 + 10u^3 = 0$.		\\
\label{Random_RS}

\indent \hspace{-5mm} Rank 2 of the Picard group:
\item $x^3 + y^3 + 2z^3 + 4u^3 = 0$.		
\label{S1}
\item $x^3 + y^3 + 5z^3 + 25u^3 = 0$.		
\label{S2}
\item $x^3 + y^3 + 3z^3 + 9u^3 = 0$.		
\label{S3}

\vspace{15mm}
\indent \hspace{-5mm}Rank 3 of the Picard group:
\item $x^3 + y^3 + 2z^3 + 2u^3 = 0$.		
\label{S4}
\item $x^3 + y^3 + 5z^3 + 5u^3 = 0$.		
\label{S5}
\item $x^3 + y^3 + 7z^3 + 7u^3 = 0$.		
\label{S6}
\item $2x^3 + 2y^3 + 3z^3 + 3u^3 = 0$.		
\label{S7}

\vspace{7mm}
\end{enumerate}
 
The first six cubic surfaces illustrate the case of Picard group rank $1$. The third surface was
considered as an example of a diagonal cubic surface with bigger coefficients. The lack of success
in finding a generating set for this surface (as opposed to all the other surfaces examined 
by that point) motivated the study of the surfaces 4--5, which have coefficients of 
intermediate value between the coefficients of the first, successful surface, and the third,
problematic one. Surface \ref{Random_RS} is aimed to illustrate the case of surfaces with 'random' 
coefficients. The remaining seven surfaces were taken from \cite{P&T} as examples of cubic surfaces with the 
rank of the Picard Group $2$ and $3$. 

In order to find a suitable generating set $G$ to begin with, we tested several small sets 
for finite generation up to a small index $n$ ($n = 100$, or $n=1000$).
We observed that, if the set $G$ generates more than $80\%-90\%$ of the first $n$ points
for a small $n$,  
then this is a good indicator that the set $G$ will generate roughly the same percentage of all
points up to a much bigger index bound $N$ (which we took to be either $5\cdot 10^4$ or $10^5$).
We chose the initial small sets to be the set of points of indexes $\{1,2,3,4\}$. If this did not
yield a large enough percentage of points generated, we would enlarge the initial set to $G = \{1,2,3,4,5\}$,
and continue this way. 
Generally, we were 'lucky', in the sense that a few tries would provide us with a good generating set $G$
(a set $G$ that generates most of the first $n$ points.) Then we would eliminate from $G$ the 
'superfluous' points, i.e., the points that could be obtained by composing other points in $G$. This is the
reason for which, for example, the first surface has $G = \{3\}$ instead of $G = \{1,2,3,4\}$: the points 
of indices 1, 2 and 4 lie in the tangent plane at the point of index $3$.
   
At first, the only
exception was the surface \ref{Big_RS}, which represents, at least computationally, a problem. 
Having added the surfaces 4--5, we noticed that it is hard to find a generating set using this
naive method for these surfaces as well.     

We found the following generating sets, listed both as sets of indices and as sets of
rational points. Here, and in all subsequent tables, the label '{\em S}' stands for 'surface'. 
\[ \begin{array}{|c||l|l|} \hline
        S & G \text{ as set of indices} & G \text { as set of points} \\
        \hline \hline
        1  &  \{3\}         &  \{ \, (1:-1:-1:1) \, \} \\
        2  & \{1,2,4\}      &  \{ \, (0:1:1:-1), (1:1:-1:0), (2:-2:1:1) \, \} \\      
        6  &  \{2\}         &  \{ \, (1:-1:-1:1) \, \} \\
       \hline
        7  &  \{3\}         &  \{ \, (1:-1:-1:1) \, \} \\     
        8  &  \{1,2\}       &  \{ \, (1:-1:0:0), (1:4:-2:-1) \, \} \\       
        9  & \{1,2,4\}      &  \{ \, (1:-1:0:0), (1:2:0:-1), (1:2:-3:2) \, \} \\
       \hline
        10  &  \{5,6\}      &  \{ \, (1:-1:-1:1), (1:-1:1:-1) \, \} \\
        11  &  \{3,4\}      &  \{ \, (1:-1:-1:1), (1:-1:1:-1) \, \} \\
        12  & \{1,2,5,6\}   &  \{ \, (0:0:1:-1), (1:-1:0:0), (1:-1:-1:1), (1:-1:1:-1) \, \} \\
        13  & \{1,2,3,4,5\} &  \{ \, (0:0:1:-1), (1:-1:0:0), (1:-1:-1:1),  \\
            &               &  (1:-1:1:-1), (3:-6:1:5) \, \} \\
       \hline
     \end{array}
\]
 
Before we go on and list the results we obtained using the {\em TestGeneratingSet} procedure, 
let us provide an indication of the asymptotics of the number of points on each cubic surface 
up to some height $H.$ Note that, as we used Dan Bernstein's program to find rational points on
the diagonal cubic surfaces, here 'height' refers to $h_{max}.$ 
The asymptotics of the number of points seems to be related with the 
percentage of points that can be generated up to some height. 
For the last seven surfaces, we did not take into consideration the 
points on the trivial rational lines, i.e., points of the type $(x:-x:y:-y)$, except for the point
$(1:-1:0:0)$ on the surfaces 7--9 and the points $(1:-1:0:0)$, $(0:0:1:-1)$, $(1:-1:1:-1)$ 
and $(1:-1:-1:1)$ on the surfaces 10--13, which we need for finite generation.

We include intermediate results of the number of points up to different height limits. These
results seem to confirm Manin's conjecture relating the 
asymptotics of rational points of bounded height on a Fano variety with the rank of the
Picard group of the surface (see \cite{FMT}:)
\[ \# \{P \in V(\Q) : h(P) < H \} \,\, \sim  \,\, CH \log ^{rkPic(V) - 1}H \]  
for $H \longrightarrow \infty$, where $h$ is an anticanonical height on $V$.

\[ \begin{array}{|c||r|r|r|r|r|r|r|r|r|r|} \hline
        S & \multicolumn{10}{|c|}{ \text {Number of points up to height} }
        \\
        & 100 & 200 & 500 & 1000 & 2000 & 5000
        		& 10000 & 20000 & 50000 & 100000 
        \\   \hline \hline
        1  &  77 &  163 & 436 & 906 & 1827	& 4408 & 8754 & 17332 & 43280 & 86329  \\
        2  & 180 &	358 &	855	& 1683 & 3244	& 8097 & 16436 & 32704 & 82581 & 166825 \\
				3  & 16	 &  25	& 62	& 117	 & 204	& 502	 & 1055	 & 2084	 & 5479 &	10840 \\
				4  & 37  &	78	& 206	& 414  & 778  &	1937 & 3877  & 7756 &	19701 & 39433 \\
				5  & 37	 &  67	& 165	& 310	 & 595	& 1580 & 3148	& 6257 & 15499 & 31134 \\
				6  & 55  & 	120 &	316 &	646  & 1285 &	3131 & 6397 &	12753 &	32072 &	64102	\\
				\hline			
				7  & 196 &	458 &	1308 & 2746 & 6004 & 16758 & 35958 & 75984 & 205284 & 433526 \\
				8  & 142 &	292 &	766 &	1734 & 3872 &	10892 &	23338 &	49608 &	135128 & 286040 \\
				9  & 200 &	438 &	1270 & 2768 & 6200 & 17434 & 37018 & 78980 & 215626 &	455164 \\
				\hline
				10  & 666 &	1630 & 5410 & 12870 & 29926 &	89218 &	205198 & 465226 &	1364810 & 3051198 \\
				11 & 412 &	1012 & 3328 &	7964 & 18676 & 56412 & 131512 &	299776 & 881774 &	1976482 \\
				12  & 702 &	1870 & 6010 &	14130 &	33156 &	100580 & 228696	& 520700 & 1526532 & 3420784 \\
				13 & 384 &	1052 & 3196 &	7752 & 18400 & 56348 & 130476 &	298860 & 876776 &	1966160 \\
       \hline
     \end{array}
\]

For the surfaces with rank of the Picard group equal to $1$ we computed, additionally,
the number of rational points up to slightly greater height limits, as summarized in 
the table below ('-' means 'not computed'.)
   
\[ \begin{array}{|r||r|r|r|r|} \hline
        S & \multicolumn{4}{|c|}{ \text {Number of points up to height} }
        \\
        & 150000 & 200000 & 250000 & 300000 
        \\   \hline \hline
				1  & 129473 & - & - & - \\ 
 				2  & 250286 & - & - & -\\
				3  & 16123 & 21627 & 27026 & 32507 \\
				4  & 59100 & 78498 & - & - \\
				5  & 46436 & 61958 & 77518 & 93079 \\
				6  & 96065 & - & - & - \\
				\hline
   \end{array}
\]

Relevant to our claim that these results seem to confirm Manin's conjecture are the 
following graphs based on the tables above. In all graphs, we plotted the number of points
up to height $H$ divided by $H \log ^{rkPic(V) - 1}H$ for various values of $H.$ The conjecture would
be verified if the plotted points would become arbitrarily close, in the limit, to a line
parallel to the $Ox$ axis, of equation $y=C,$ where $C$ is the constant predicted by 
Manin's conjecture. For a conjecture about the value of this constant, see \cite{P&T}. 

\begin{figure}
	\centering
		\includegraphics[width=1.21\textwidth]{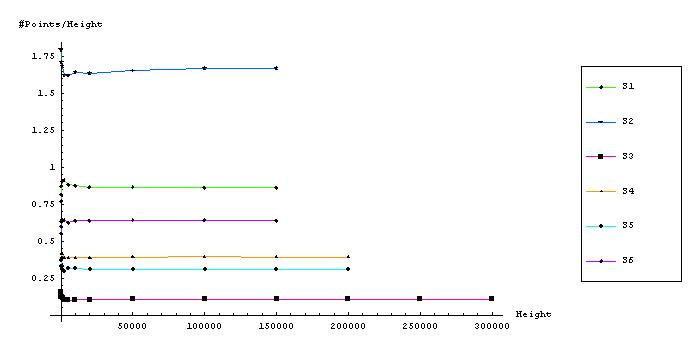}
		\caption{Surfaces with Picard group rank $1$}
	\label{r1}
\end{figure}

\begin{figure}
	\centering
		\includegraphics[width=1.32\textwidth]{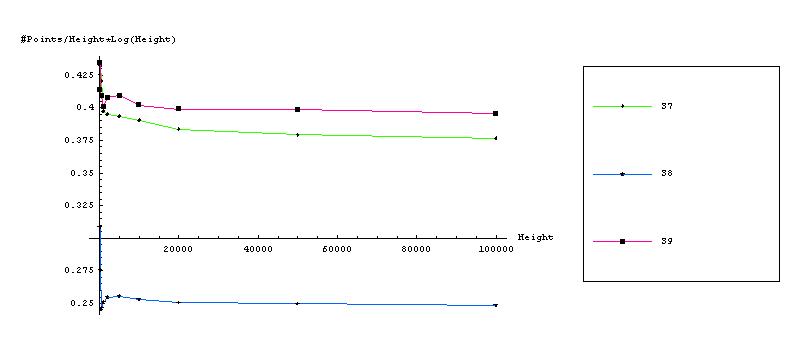}
		\caption{Surfaces with Picard group rank $2$}
	\label{r2}
\end{figure}

\begin{figure}
	\centering
		\includegraphics[width=1.27\textwidth]{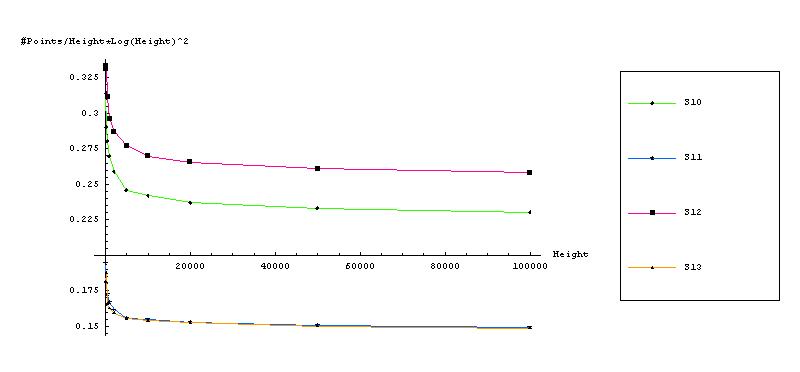}
		\caption{Surfaces with Picard group rank $3$}
	\label{r3}
\end{figure}

In the remaining, by 'height' we mean $h_{sum}$.

Note that for the surfaces with rank of the Picard group equal to two, 
most of the points are 'doubled', i.e., if $(x:y:z:u)$ is a point
on the cubic surface, then so is $(y:x:z:u)$, while for the surfaces
with rank of the Picard group equal to three, most of the points are 'quadrupled',
i.e., if $(x:y:z:u)$ is a point
on the cubic surface, then so are $(y:x:z:u)$, $(x:y:u:z)$ and $(y:x:u:z)$. 
In the following we list the results which were obtained using the {\em TestGeneratingSet} 
procedure. The generating sets used are the ones enumerated above, while the index bounds 
and the corresponding height bounds are 
given in the third and second columns of the table. '$\#$ iter' is the number of iterations
of the outer loop of the procedure, and the 'first bad point' refers to the point of 
smallest index that could not be generated by the procedure. For example, the first line
in the table reads ''The procedure {\em TestGeneratingSet} called for surface 1, with
index bound $100$ corresponding to the height bound $317$,
and initial generating set $G = \{3\}$ (or $G = \{ \, (1:-1:-1:1) \, \}$),
generates $74$ rational points, which represents 
$74.0\%$ of the first $100$ points, in $4$ iterations
of the outer loop. The smallest point that could not be 
generated has index $30$ and height $86$.'' 
   
\[ \begin{array}{|c||r|r|r|c|c|r|r|r|} \hline
       \text{Surface} & \text{Height} & \text{Index} &  \text{$\#$ points} & \text{$\%$ points} 
       			& \text{$\#$ iter} & \multicolumn{2}{|c|} { \text{First bad point} }
        \\
        & \text{bound} & \text{bound} & \text{generated} & \text{generated} & & \text{Index} & \text{Height}
        \\   \hline \hline
        1  & 317 & 100  & 74 & 74.0 & 4 & 30 & 86  \\
        1  & 617 & 200 & 160 & 80.0 & 9 & 30 & 86  \\
        1  & 1,443 & 500 & 463 & 92.6 & 16 & 42 & 130  \\
        1  & 2,788 & 1,000 & 923 & 92.3 & 15 & 255 & 788  \\
        1  & 5,574 & 2,000  & 1,859 & 93.0 & 14 & 543 & 1,541  \\
        1  & 14,456 & 5,000  & 4,747 & 94.9 & 15 & 1,145 & 3,192  \\
        1  & 29,074 & 10,000 & 9,462 & 94.6 & 14 & 1,593 & 4,423  \\
        1  & 58,775 & 20,000 & 18,957 & 94.8 & 14 & 3,633 & 10,322  \\
        1  & 147,343 & 50,000  & 47,418 & 94.8 & 13 & 8,522 & 24,677  \\
        1  & 296,822 & 100,000 & 94,910 & 94.9 & 13 & 8,522 & 24,677  \\
        \hline
        2  & 150 & 100 & 97 & 97.0 & 7 & 85 & 124  \\
        2  & 282 & 200 & 196 & 98.0 & 9 & 90 & 134  \\
        2  & 703 & 500 & 483 & 96.6 & 8 & 258 & 364  \\
        2  & 1,477 & 1,000 & 973 & 97.3 & 9 & 358 & 511  \\
        2  & 3,020 & 2,000 & 1931 & 96.6 & 9 & 625 & 943  \\
        2  & 7,663 & 5,000 & 4,813 & 96.3 & 10 & 1,040 & 1,542  \\
        2  & 15,405 & 10,000 & 9,659 & 96.6 & 11 & 1,775 & 2,656  \\
        2  & 30,651 & 20,000 & 19,259 & 96.3 & 11 & 4,262 & 6,539  \\
        2  & 75,845 & 50,000 & 48,181 & 96.3 & 11 & 10,073 & 15,539  \\
        2  & 151,171 & 100,000 & 96,477 & 96.5 & 12 & 15,223 & 23,243  \\ 
        \hline       
        6  & 388 & 100 & 86 & 86.0 & 5 & 49 & 209  \\
        6  & 762 & 200 & 176 & 88.0 & 5 & 49 & 209  \\
        6  & 1,864 & 500 & 468 & 93.6 & 10 & 169 & 641  \\
        6  & 3,687 & 1,000 & 937 & 93.7 & 11 & 181 & 688  \\
        6  & 7,557 & 2,000 & 1,867 & 93.3 & 11 & 513 & 1,926  \\
        6  & 18,976 & 5,000 & 4,677 & 93.6 & 11 & 1,078 & 3,984  \\
        6  & 37,612 & 10,000 & 9,410 & 94.1 & 11 & 2,271 & 8,661  \\
        6  & 74,617 & 20,000 & 18,963 & 94.8 & 11 & 2,662 & 10,125  \\
        6  & 186,532 & 50,000 & 47,436 & 94.9 & 12 & 6,373 & 24,068  \\
        \hline       
     \end{array}
\]

\[ \begin{array}{|c||r|r|r|c|c|r|r|r|} \hline
       \text{Surface} & \text{Height} & \text{Index} &  \text{$\#$ points} & \text{$\%$ points} 
       			& \text{$\#$ iter} & \multicolumn{2}{|c|} { \text{First bad point} }
        \\
        & \text{bound} & \text{bound} & \text{generated} & \text{generated} & & \text{Index} & \text{Height}
        \\   \hline \hline
        7  & 129 & 100 & 100 & 100.0 & 6 & - & -  \\
        7  & 245 & 200 & 194 & 97.0 & 6 & 127 & 167  \\
        7  & 538 & 500 & 490 & 98.0 & 8 & 304 & 376  \\
        7  & 980 & 1,000 & 990 & 99.0 & 7 & 550 & 612  \\
        7  & 1,889 & 2,000 & 1,984 & 99.2 & 7 & 1,022 & 992  \\
        7  & 4,230 & 5,000 & 4,974 & 99.5 & 7 & 2,620 & 2,401  \\
        7  & 7,974 & 10,000 & 9,934 & 99.3 & 8 & 5,610 & 4,707  \\
        7  & 14,775 & 20,000 & 19,934 & 99.7 & 8 & 7,512 & 6,222 \\
        7  & 34,339 & 50,000 & 49,880 & 99.8 & 7 & 19,666 & 14,554 \\
        7  & 64,682 & 100,000 & 99,812 & 99.8 & 8 & 38,212 & 26,672  \\
        7  & 94,215 & 150,000 & 149,744 & 99.8 & 9 & 38,212 & 26,672  \\
        \hline       
        8  & 172 & 100 & 81 & 81.0 & 6 & 42 & 78  \\
        8  & 316 & 200 & 170 & 85.0 & 8 & 56 & 104  \\
        8  & 750 & 500 & 488 & 97.6 & 9 & 152 & 234  \\
        8  & 1,412 & 1,000 & 988 & 98.8 & 8 & 516 & 774  \\
        8  & 2,484 & 2,000 & 1,960 & 98.0 & 8 & 516 & 774  \\
        8  & 5,632 & 5,000 & 4,922 & 98.4 & 9 & 1,855 & 2,322  \\
        8  & 10,354 & 10,000 & 9,874 & 98.7 & 8 & 3,708 & 4,296  \\
        8  & 19,444 & 20,000 & 19,836 & 99.2 & 8 & 6,852 & 7,538  \\
        8  & 44,750 & 50,000 & 49,720 & 99.4 & 8 & 16,058 & 15,812  \\
        8  & 84,436 & 100,000 & 99,626 & 99.6 & 9 & 32,420 & 30,072  \\
        \hline
        9  & 114 & 100 & 48 & 48.0 & 4 & 8 & 24  \\
        9  & 242 & 200 & 198 & 99.0 & 9 & 126 & 146  \\
        9  & 522 & 500 & 484 & 96.8 & 9 & 318 & 346  \\
        9  & 978 & 1,000 & 956 & 95.6 & 11 & 379 & 414  \\
        9  & 1,822 & 2,000 & 1,968 & 98.4 & 9 & 781 & 770  \\
        9  & 3,878 & 5,000 & 4,954 & 99.1 & 9 & 1,602 & 1,472  \\
        9  & 7,254 & 10,000 & 9,936 & 99.4 & 9 & 3,728 & 3,046  \\
        9  & 13,610 & 20,000 & 19,908 & 99.5 & 9 & 10,420 & 7,522  \\
        9  & 31,320 & 50,000 & 49,806 & 99.6 & 8 & 21,142 & 14,342  \\
        9  & 58,852 & 100,000 & 99,778 & 99.8 & 9 & 32,036 & 20,884  \\
        \hline \hline
        10  & 61 & 100 & 92 & 92.0 & 3 & 79 & 51  \\
        10  & 91 & 200 & 200 & 100.0 & 3 & - & -  \\
        10  & 214 & 500 & 496 & 99.2 & 5 & 419 & 184  \\
        10  & 358 & 1,000 & 980 & 98.0 & 6 & 651 & 255  \\
        10  & 612 & 2,000 & 1,996 & 99.8 & 5 & 1,791 & 554  \\
        10  & 1,225 & 5,000 & 4,940 & 98.8 & 5 & 2,259 & 674  \\
        10  & 2,143 & 10,000 & 9,916 & 99.2 & 6 & 3,675 & 976  \\
        10  & 3,806 & 20,000 & 19,852 & 99.3 & 6 & 5,779 & 1,396  \\
        10  & 8,020 & 50,000 & 49,732 & 99.5 & 7 & 20,870 & 3,949  \\
        \hline        
     \end{array}
\]

\[ \begin{array}{|c||r|r|r|c|c|r|r|r|} \hline
       \text{Surface} & \text{Height} & \text{Index} &  \text{$\#$ points} & \text{$\%$ points} 
       			& \text{$\#$ iter} & \multicolumn{2}{|c|} { \text{First bad point} }
        \\
        & \text{bound} & \text{bound} & \text{generated} & \text{generated} & & \text{Index} & \text{Height}
        \\   \hline \hline
        11  & 94 & 100 & 89 & 89.0 & 5 & 61 & 56  \\
        11  & 144 & 200 & 184 & 92.0 & 5 & 61 & 56  \\
        11  & 274 & 500 & 492 & 98.4 & 6 & 257 & 174  \\
        11  & 474 & 1,000 & 988 & 98.8 & 6 & 757 & 382  \\
        11  & 802 & 2,000 & 1,960 & 98.0 & 6 & 1,177 & 528  \\
        11  & 1,688 & 5,000 & 4,924 & 98.5 & 7 & 1,495 & 642  \\
        11  & 2,882 & 10,000 & 9,888 & 98.9 & 8 & 3,873 & 1,386  \\
        11  & 5,100 & 20,000 & 19,732 & 98.7 & 9 & 6,207 & 2,004  \\
        11  & 10,880 & 50,000 & 49,544 & 99.1 & 9 & 11,737 & 3,308  \\
        \hline
        12  & 48 & 100 & 96 & 96.0 & 4 & 95 & 46  \\
        12  & 92 & 200 & 192 & 96.0 & 4 & 95 & 46  \\
        12  & 186 & 500 & 476 & 95.2 & 5 & 223 & 106  \\
        12  & 286 & 1,000 & 956 & 95.6 & 6 & 223 & 106  \\
        12  & 484 & 2,000 & 1,969 & 98.5 & 5 & 964 & 284  \\
        12  & 1,014 & 5,000 & 4,911 & 98.2 & 6 & 2,315 & 548  \\
        12  & 1,740 & 10,000 & 9,880 & 98.8 & 6 & 3,486 & 764  \\
        12  & 3,066 & 20,000 & 19,832 & 99.2 & 7 & 4,030 & 856  \\
        12  & 6,514 & 50,000 & 49,532 & 99.1 & 7 & 16,064 & 2,578  \\
        \hline
        13  & 106 & 100 & 96 & 96.0 & 4 & 41 & 75  \\
        13  & 167 & 200 & 196 & 98.0 & 5 & 169 & 153  \\
        13  & 316 & 500 & 484 & 96.8 & 6 & 169 & 153  \\
        13  & 515 & 1,000 & 980 & 98.0 & 6 & 572 & 360  \\
        13  & 910 & 2,000 & 1,944 & 97.2 & 6 & 860 & 465  \\
        13  & 1,885 & 5,000 & 4,896 & 97.9 & 7 & 1,937 & 897  \\
        13  & 3,310 & 10,000 & 9,780 & 97.8 & 7 & 3,102 & 1,323  \\
        13  & 5,727 & 20,000 & 19,672 & 98.4 & 7 & 4,785 & 1,816  \\       
        13  & 12,139 & 50,000 & 48,256 & 96.5 & 8 & 8,202 & 2,805  \\
        \hline
     \end{array}
\]

Note that, in general, when using 
a greater index bound we found that the 'first bad point' changed (i.e., another point
of greater height and index became the 'first bad point'), meaning that
using stepping stones of bigger height typically fills up the gaps obtained when
using a lower index bound. This is a good indicator that if we continue increasing 
the index (and thus the height) bounds, we will gradually generate {\em all} the points up to bigger and 
bigger heights. 

Let us see now what happens with the 'problematic surfaces' 3--5. Unfortunately,
any try of finding a generating set to begin with, that finds 'first bad points'
of increasing height, and that generates a percentage
of points similar to the ones obtained for the 'good' surfaces was not successful. 
Not even a 'brute force' approach like considering the initial {\em GeneratedSet}
to be, say, the first $100$ or $1000$ points does not yield satisfactory results. 
The results are better for the surfaces 4--5 than for the surface 3, with the
biggest coefficients, but still very 'bad'. Here is an illustration of the behavior
of these surfaces when starting with the {\em GeneratedSet} $= \{1,2,\ldots,10\}$:       

\[ \begin{array}{|c||r|r|r|c|c|r|r|r|} \hline
       \text{Surface} & \text{Height} & \text{Index} &  \text{$\#$ points} & \text{$\%$ points} 
       			& \text{$\#$ iter} & \multicolumn{2}{|c|} { \text{First bad point} }
        \\
        & \text{bound} & \text{bound} & \text{generated} & \text{generated} & & \text{Index} & \text{Height}
        \\   \hline \hline
        3  & 2,161 & 100 & 17 & 17.0 & 2 & 13 & 203  \\
        3  & 5,495 & 200 & 24 & 12.0 & 2 & 13 & 203  \\
        3  & 13,429 & 500 & 35 & 7.0 & 2 & 13 & 203  \\
        3  & 25,874 & 1,000 & 49 & 4.9 & 2 & 13 & 203  \\
        3  & 51,663 & 2,000 & 81 & 4.1 & 2 & 13 & 203  \\
        3  & 124,062 & 5,000 & 154 & 3.1 & 2 & 13 & 203  \\
        3  & 251,103 & 10,000 & 274 & 2.7 & 2 & 13 & 203  \\
        3  & 505,619 & 20,000 & 429 & 2.1 & 2 & 13 & 203  \\      
        \hline
        4  & 658 & 100 & 26 & 26.0 & 1 & 12 & 50  \\
        4  & 1,345 & 200 & 50 & 25.0 & 2 & 12 & 50  \\
        4  & 3,307 & 500 & 102 & 20.4 & 2 & 12 & 50  \\
        4  & 6,774 & 1,000 & 172 & 17.2 & 3 & 12 & 50  \\
        4  & 13,772 & 2,000 & 284 & 14.2 & 3 & 12 & 50  \\
        4  & 34,552 & 5,000 & 487 & 9.7 & 3 & 12 & 50  \\
        4  & 68,425 & 10,000 & 781 & 7.8 & 3 & 12 & 50  \\
        4  & 135,691 & 20,000 & 1,222 & 6.1 & 4 & 12 & 50  \\             
        \hline
        5  & 844 & 100 & 19 & 19.0 & 1 & 13 & 103  \\
        5  & 1,691 & 200 & 26 & 13.0 & 2 & 13 & 103  \\
        5  & 4,394 & 500 & 51 & 10.2 & 2 & 13 & 103  \\
        5  & 8,780 & 1,000 & 80 & 8.0 & 2 & 13 & 103  \\
        5  & 16,962 & 2,000 & 119 & 6.0 & 2 & 13 & 103  \\
        5  & 43,224 & 5,000 & 216 & 4.3 & 2 & 13 & 103  \\
        5  & 87,176 & 10,000 & 338 & 3.4 & 3 & 13 & 103  \\
        5  & 174,128 & 20,000 & 538 & 2.7 & 3 & 13 & 103  \\   
        \hline
     \end{array}
\]

These results seem to support either that $\{1,2,\ldots,10\}$ is not a generating
set for any of the three surfaces, or that the stepping stones needed to fill up the gaps
(i.e., the rational points needed to decompose the 'first bad points') have very big heights.    
Although the percentages of generated points 
obtained for the surfaces 4--5 are slightly better than the percentages for the
surface 3, they still become smaller and smaller as the index bound limit (and so
also the height) grow. But the most important negative indicator is that 'the first bad point'
never changes.    

In order to make progress, we introduced another approach to finding a generating set for 
the surfaces 3--5, based on the idea of 'throwing in' (adding to the {\em Generated Set})
the first bad points if they cannot be generated by decomposition.  
Our aim is to obtain, after adding sufficiently many 'first bad points', a set of points
that generates a stable (or even better, increasing) percentage of 
points for increasing index bounds, and a 'changing first bad point' 
behavior, i.e., applying the {\em TGS} procedure to increasing index bounds would result in 
finding 'first bad points' of increasing heights.

We implement this new approach in the following way. We apply the {\em TGS} procedure
to a (small) generating set and an index bound of $1000$. We obtain a 'first bad point'
that unfortunately stays the same when increasing the index bound (as observed when using
our first approach). We apply again the {\em TGS} procedure
to the initial generating set {\em and} this first bad point, with an index bound of $1000$.
We obtain another 'first bad point', of bigger index and height than the initial one. We add this point
to our generating set (which now contains also the initial 'first bad point') and continue
this way. We stop when we have added sufficiently many 'first bad points' to our initial set 
so that this new, bigger generating set fulfills the two objectives mentioned above. Once we
have obtained such a set, we stop adding points to our generating set and just increase the
index bounds to make sure the percentage of generated points is indeed stable or increasing,
and that the height of the 'first bad point' grows as the index bound is increased.       
  
For example, for surface 4, we start with {\em Generated Set} $= \{1, 2, \ldots, 10 \}$. 
We obtain the first bad point $12$, which is stable - stays the same even if we increase the index bound. 
We add it to the {\em Generated Set} and call again the {\em TGS} procedure. We obtain more points, and another 
first bad point. We add this new bad point to the {\em Generated Set} and continue this way, gradually
filling the holes. At first we kept the index bound constant, until we obtained a reasonable percentage
of generated points. Then we tested whether the 'first bad point' changes when increasing the index 
bound and keeping the initial {\em Generated Set} constant (i.e., we stopped filling the holes, and
just increased the index bound.) For surfaces 4 and 5 this approach seems successful, as reflected in the
tables below.

\[ \begin{array}{|c||r|r|r|c|c|r|r|r|} \hline
       \text{Surface} & \text{Height} & \text{Index} &  \text{$\#$ points} & \text{$\%$ points} 
       			& \text{$\#$ iter} & \multicolumn{2}{|c|} { \text{First bad point} }
        \\
        & \text{bound} & \text{bound} & \text{generated} & \text{generated} & & \text{Index} & \text{Height}
        \\   \hline \hline
        4  & 6,774 & 1,000 & 172 & 17.2 & 3 & 12 & 50  \\
        4  & 6,774 & 1,000 & 177 & 17.7 & 3 & 13 & 55  \\
        4  & 6,774 & 1,000 & 194 & 19.4 & 4 & 14 & 63  \\
        4  & 6,774 & 1,000 & 210 & 21.0 & 4 & 15 & 73 \\
        4  & 6,774 & 1,000 & 218 & 21.8 & 4 & 20 & 107  \\
        4  & 6,774 & 1,000 & 230 & 23.0 & 4 & 21 & 108  \\
        4  & 6,774 & 1,000 & 237 & 23.7 & 4 & 22 & 110  \\
        4  & 6,774 & 1,000 & 249 & 24.9 & 5 & 23 & 125  \\
        4  & 6,774 & 1,000 & 268 & 26.8 & 6 & 25 & 179  \\ 
        4  & 6,774 & 1,000 & 282 & 28.2 & 6 & 27 & 193  \\
        4  & 6,774 & 1,000 & 296 & 29.6 & 6 & 28 & 199  \\
        4  & 6,774 & 1,000 & 325 & 32.5 & 13 & 32 & 215  \\
        4  & 6,774 & 1,000 & 328 & 32.8 & 13 & 35 & 249  \\
        4  & 6,774 & 1,000 & 335 & 33.5 & 13 & 37 & 262  \\
        4  & 6,774 & 1,000 & 338 & 33.8 & 13 & 43 & 297  \\
        4  & 6,774 & 1,000 & 342 & 34.2 & 13 & 49 & 317  \\ 
        4  & 6,774 & 1,000 & 349 & 34.9 & 13 & 52 & 329  \\
        4  & 6,774 & 1,000 & 351 & 35.1 & 13 & 58 & 370 \\
        4  & 6,774 & 1,000 & 353 & 35.3 & 13 & 62 & 396  \\
        4  & 6,774 & 1,000 & 360 & 36.0 & 13 & 66 & 413  \\
        4  & 6,774 & 1,000 & 372 & 37.2 & 13 & 69 & 438  \\
        4  & 6,774 & 1,000 & 394 & 39.4 & 18 & 73 & 467  \\
        4  & 6,774 & 1,000 & 400 & 40.0 & 18 & 76 & 487  \\              
        \hline
        4  & 34,552 & 5,000 & 1,331 & 26.6 & 38 & 89 & 570  \\
        4  & 68,425 & 10,000 & 2,769 & 27.7 & 50 & 92 & 611  \\
        4  & 135,691 & 20,000 & 6,365 & 31.8 & 53 & 189 & 1,230  \\
        4  & 204,042 & 30,000 & 10,142 & 33.8 & 50 & 233 & 1,605  \\
        4  & 271,092 & 40,000 & 14,403 & 36.0 & 45 & 324 & 2,115  \\
        4  & 339,994 & 50,000 & 18,409 & 36.8 & 51 & 352 & 2,387  \\
        \hline
        \end{array}
\] 
        
\[ \begin{array}{|c||r|r|r|c|c|r|r|r|} \hline
       \text{Surface} & \text{Height} & \text{Index} &  \text{$\#$ points} & \text{$\%$ points} 
       			& \text{$\#$ iter} & \multicolumn{2}{|c|} { \text{First bad point} }
        \\
        & \text{bound} & \text{bound} & \text{generated} & \text{generated} & & \text{Index} & \text{Height}
        \\   \hline \hline
        5  & 8,780 & 1,000 & 80 & 8.0 & 2 & 13 & 103  \\
        5  & 8,780 & 1,000 & 87 & 8.7 & 3 & 14 & 111  \\
        5  & 8,780 & 1,000 & 100 & 10.0 & 3 & 15 & 112  \\
        5  & 8,780 & 1,000 & 114 & 11.4 & 6 & 16 & 122 \\
        5  & 8,780 & 1,000 & 142 & 14.2 & 8 & 17 & 125  \\
        5  & 8,780 & 1,000 & 149 & 14.9 & 8 & 18 & 126  \\
        5  & 8,780 & 1,000 & 157 & 15.7 & 8 & 19 & 127  \\ 
        5  & 8,780 & 1,000 & 170 & 17.0 & 8 & 21 & 150 \\
        5  & 8,780 & 1,000 & 175 & 17.5 & 8 & 23 & 168  \\
        5  & 8,780 & 1,000 & 177 & 17.7 & 8 & 25 & 177  \\
        5  & 8,780 & 1,000 & 207 & 20.7 & 16 & 27 & 188  \\         
        5  & 8,780 & 1,000 & 211 & 21.1 & 16 & 28 & 190  \\ 
        5  & 8,780 & 1,000 & 219 & 21.9 & 16 & 32 & 211  \\       
        5  & 8,780 & 1,000 & 223 & 22.3 & 16 & 37 & 276  \\  
        5  & 8,780 & 1,000 & 230 & 23.0 & 16 & 39 & 298  \\
        5  & 8,780 & 1,000 & 232 & 23.2 & 16 & 44 & 350  \\         
        5  & 8,780 & 1,000 & 236 & 23.6 & 16 & 45 & 363  \\
        5  & 8,780 & 1,000 & 237 & 23.7 & 16 & 46 & 367  \\ 
        5  & 8,780 & 1,000 & 242 & 24.2 & 16 & 47 & 369  \\         
        5  & 8,780 & 1,000 & 268 & 26.8 & 16 & 56 & 427  \\        
        5  & 8,780 & 1,000 & 276 & 27.6 & 16 & 57 & 431  \\        
        5  & 8,780 & 1,000 & 282 & 28.2 & 16 & 59 & 445  \\        
        5  & 8,780 & 1,000 & 311 & 31.1 & 16 & 60 & 464  \\        
        5  & 8,780 & 1,000 & 313 & 31.3 & 16 & 62 & 487  \\         
        5  & 8,780 & 1,000 & 319 & 31.9 & 16 & 66 & 581  \\               
        5  & 8,780 & 1,000 & 337 & 33.7 & 16 & 68 & 595  \\        
        5  & 8,780 & 1,000 & 339 & 33.9 & 16 & 69 & 602  \\        
        5  & 8,780 & 1,000 & 347 & 34.7 & 16 & 75 & 631  \\                
        5  & 8,780 & 1,000 & 356 & 35.6 & 16 & 76 & 637  \\                       
        5  & 8,780 & 1,000 & 365 & 36.5 & 16 & 84 & 695  \\
        5  & 8,780 & 1,000 & 369 & 36.9 & 16 & 87 & 719  \\               
        5  & 8,780 & 1,000 & 380 & 38.0 & 16 & 88 & 733  \\                       
        5  & 8,780 & 1,000 & 385 & 38.5 & 16 & 91 & 745  \\        
        5  & 8,780 & 1,000 & 390 & 39.0 & 16 & 93 & 771  \\        
        5  & 8,780 & 1,000 & 409 & 40.9 & 16 & 96 & 801  \\                      
        5  & 8,780 & 1,000 & 413 & 41.3 & 16 & 103 & 862  \\
        \hline
        5  & 43,224 & 5,000 & 1,881 & 37.6 & 29 & 118 & 1,015  \\
        5  & 87,176 & 10,000 & 3,650 & 36.5 & 32 & 145 & 1,197  \\ 
        5  & 174,128 & 20,000 & 7,236 & 36.2 & 37 & 295 & 2,554  \\
        5  & 262,052 & 30,000 & 11,367 & 37.9 & 44 & 325 & 2,774  \\
        5  & 349,121 & 40,000 & 15,842 & 39.6 & 44 & 461 & 3,988  \\
        5  & 437,046 & 50,000 & 20,103 & 40.2 & 35 & 461 & 3,988  \\
        \hline
     \end{array}
\] 

Unfortunately, for surface 3 this approach does not seem to work. After adding 
many more 'first bad points' to the initial generating set than for the 
surfaces 4--5, we still did not obtain a 'good' generating set, as illustrated
below. 

\[ \begin{array}{|c||r|r|r|c|c|r|r|r|} \hline
       \text{Surface} & \text{Height} & \text{Index} &  \text{$\#$ points} & \text{$\%$ points} 
       			& \text{$\#$ iter} & \multicolumn{2}{|c|} { \text{First bad point} }
        \\
        & \text{bound} & \text{bound} & \text{generated} & \text{generated} & & \text{Index} & \text{Height}
        \\   \hline \hline
        3  & 25,874 & 1,000 & 49 & 4.9 & 2 & 13 & 203  \\
        3  & 25,874 & 1,000 & 52 & 5.2 & 2 & 14 & 248  \\
        3  & 25,874 & 1,000 & 54 & 5.4 & 2 & 15 & 260  \\        
        3  & 25,874 & 1,000 & 57 & 5.7 & 2 & 16 & 264  \\        
        3  & 25,874 & 1,000 & 60 & 6.0 & 2 & 18 & 335  \\       
        3  & 25,874 & 1,000 & 62 & 6.2 & 2 & 19 & 337  \\        
        3  & 25,874 & 1,000 & 64 & 6.4 & 2 & 20 & 383  \\       
        3  & 25,874 & 1,000 & 66 & 6.6 & 2 & 21 & 413  \\       
        3  & 25,874 & 1,000 & 67 & 6.7 & 2 & 22 & 433  \\         
        3  & 25,874 & 1,000 & 69 & 6.9 & 2 & 23 & 434  \\        
        3  & 25,874 & 1,000 & 71 & 7.1 & 2 & 26 & 526  \\       
        3  & 25,874 & 1,000 & 73 & 7.3 & 2 & 27 & 573  \\        
        3  & 25,874 & 1,000 & 76 & 7.6 & 2 & 28 & 605  \\                 
        3  & 25,874 & 1,000 & 77 & 7.7 & 2 & 29 & 630  \\                
        3  & 25,874 & 1,000 & 78 & 7.8 & 2 & 31 & 699  \\               
        3  & 25,874 & 1,000 & 80 & 8.0 & 2 & 32 & 711  \\         
        3  & 25,874 & 1,000 & 82 & 8.2 & 2 & 35 & 754  \\        
        3  & 25,874 & 1,000 & 85 & 8.5 & 2 & 36 & 772 \\        
        3  & 25,874 & 1,000 & 86 & 8.6 & 2 & 37 & 775  \\                      
        3  & 25,874 & 1,000 & 88 & 8.8 & 2 & 39 & 808  \\                 
        3  & 25,874 & 1,000 & 90 & 9.0 & 2 & 40 & 819  \\               
        3  & 25,874 & 1,000 & 93 & 9.3 & 2 & 41 & 853 \\       
        3  & 25,874 & 1,000 & 95 & 9.5 & 2 & 42 & 868  \\                
        3  & 25,874 & 1,000 & 98 & 9.8 & 2 & 43 & 872 \\               
        3  & 25,874 & 1,000 & 99 & 9.9 & 2 & 44 & 895  \\                              
        3  & 25,874 & 1,000 & 100 & 10.0 & 2 & 45 & 895  \\                        
        3  & 25,874 & 1,000 & 106 & 10.6 & 3 & 48 & 1,021  \\                       
        3  & 25,874 & 1,000 & 108 & 10.8 & 3 & 49 & 1,032 \\       
        3  & 25,874 & 1,000 & 109 & 10.9 & 3 & 50 & 1,042 \\        
        3  & 25,874 & 1,000 & 110 & 11.0 & 3 & 51 & 1,061  \\                                     
        3  & 25,874 & 1,000 & 111 & 11.1 & 3 & 52 & 1,062  \\                               
        3  & 25,874 & 1,000 & 112 & 11.2 & 3 & 53 & 1,079  \\                               
        3  & 25,874 & 1,000 & 113 & 11.3 & 3 & 54 & 1,097 \\               
        3  & 25,874 & 1,000 & 116 & 11.6 & 3 & 55 & 1,120 \\       
        3  & 25,874 & 1,000 & 117 & 11.7 & 3 & 56 & 1,131 \\               
        3  & 25,874 & 1,000 & 118 & 11.8 & 3 & 58 & 1,226  \\                                            
        3  & 25,874 & 1,000 & 120 & 12.0 & 3 & 59 & 1,270  \\                   
        \hline
        \end{array}
\] 

Since this is going way too slow, we will 'throw' in our {\em Generated Set} not only the first bad point, 
but the first 10 bad points. 

\[ \begin{array}{|c||r|r|r|c|c|r|r|r|} \hline
       \text{Surface} & \text{Height} & \text{Index} &  \text{$\#$ points} & \text{$\%$ points} 
       			& \text{$\#$ iter} & \multicolumn{2}{|c|} { \text{First bad point} }
        \\
        & \text{bound} & \text{bound} & \text{generated} & \text{generated} & & \text{Index} & \text{Height}
        \\   \hline \hline
        3  & 25,874 & 1,000 & 137 & 13.7 & 3 & 69 & 1,496  \\       
        3  & 25,874 & 1,000 & 151 & 15.1 & 3 & 82 & 1,741  \\     
        3  & 25,874 & 1,000 & 164 & 16.4 & 3 & 95 & 2,110  \\      
        3  & 25,874 & 1,000 & 177 & 17.7 & 2 & 107 & 2,458  \\        
        3  & 25,874 & 1,000 & 187 & 18.7 & 2 & 118 & 2,753  \\       
        3  & 25,874 & 1,000 & 207 & 20.7 & 5 & 134 & 3,039  \\            
        3  & 25,874 & 1,000 & 220 & 22.0 & 5 & 146 & 3,391  \\             
        3  & 25,874 & 1,000 & 233 & 23.3 & 5 & 160 & 3,928  \\      
        3  & 25,874 & 1,000 & 243 & 24.3 & 5 & 174 & 4,686  \\              
        3  & 25,874 & 1,000 & 255 & 25.5 & 5 & 184 & 4,865  \\                
        3  & 25,874 & 1,000 & 268 & 26.8 & 5 & 197 & 5,257  \\
        \hline
        \end{array}
\] 

This is again too slow, so we start inserting the first 20 bad points to our {\em Generated Set}. 

\[ \begin{array}{|c||r|r|r|c|c|r|r|r|} \hline
       \text{Surface} & \text{Height} & \text{Index} &  \text{$\#$ points} & \text{$\%$ points} 
       			& \text{$\#$ iter} & \multicolumn{2}{|c|} { \text{First bad point} }
        \\
        & \text{bound} & \text{bound} & \text{generated} & \text{generated} & & \text{Index} & \text{Height}
        \\   \hline \hline        
        3  & 25,874 & 1,000 & 301  & - & 30.1 & 5 & 226 & 6,309  \\       
        3  & 25,874 & 1,000 & 325  & - & 32.5 & 5 & 248 & 6,811  \\        
        3  & 25,874 & 1,000 & 347  & - & 34.7 & 5 & 269 & 7,255  \\         
        3  & 25,874 & 1,000 & 367  & - & 36.7 & 5 & 290 & 7,873  \\               
        3  & 25,874 & 1,000 & 388  & - & 38.8 & 5 & 314 & 8,592  \\       
        3  & 25,874 & 1,000 & 409  & - & 40.9 & 5 & 338 & 9,134  \\                     
        3  & 25,874 & 1,000 & 434  & - & 43.4 & 5 & 359 & 9,673  \\             
        \hline
        3  & 51,663 & 2,000 & 536 & 436 & 26.8 & 5 & 359 & 9,673  \\       
        3  & 124,062 & 5,000 & 734 & 437 & 14.7 & 5 & 359 & 9,673  \\        
        3  & 251,103 & 10,000 & 985 & 437 & 9.9 & 5 & 359 & 9,673  \\       
        3  & 505,619 & 20,000 & 1,298 & 437 & 6.5 & 7 & 359 & 9,673  \\
        \hline
        \end{array}
\] 


Next we present other statistical data. 

It seems that the percentage of points on a surface that can be strongly decomposed (a point $x$ is {\em strongly 
decomposable} if it has a decomposition $x = y \circ z$ with $y,z < x$) up to some index $N$ is approximately
constant for various values of $N$. This suggests that this percentage may be an invariant for the surface. 

It seems likely that if this percentage is bigger, than {\em TestGeneratingSet} will generate more points (up to
some index), using a suitable $GeneratedSet$. This is confirmed if we study the first two surfaces. Surface \ref{RS} 
has roughly $\frac{N}{8}$ points that are not strongly decomposable up to the index $N$ (for $N \geq 1000$), while
the surface \ref{Second_RS} has only $\sim \frac{N}{11}$ such points; and indeed, if we compare the 
results of $TGS$ for the two surfaces, we notice that $TGS$ for the surface \ref{Second_RS} generates more points
(up to the same index) than $TGS$ for the surface \ref{RS}. Also, note that the percentage of points that are strongly decomposable for the surface \ref{Big_RS} is very small (approximately $10\%$.) This may be one of 
the explanations for our lack of success with this surface.

\section{Conclusion}
 
The theory surrounding the Mordell-Weil problem for cubic surfaces seems not very well 
developed, mainly because of the difficulties caused by the lack of a group structure on
the operation of composing points. In this paper we presented numerical data for thirteen
diagonal cubic surfaces, in the hope of developing some intuition on a possible finiteness
conjecture
(first mentioned by Manin, cf. \cite {Man1} and \cite{Man}). For each of the surfaces, we tried to find a 
generating set. A naive method gave positive results for ten of the surfaces, while a
more rigorous method was needed to obtain similar (but not as positive) results for two 
of the other surfaces. For these surfaces, the numerical data suggest that they might be 
indeed finitely generated. The remaining surface resisted to both methods. We cannot say, 
however, whether this means that the surface is not finitely generated, or that this is 
just a sign of the limits of the methods used. 

{\bf Acknowledgement} The results of this paper arose as part of the 
author's research project at the Max Planck Institute for Mathematics in Bonn,
under the guidance of Y.~I.~Manin. 
The author would like to thank Y.~I.~Manin for providing this very enjoyable opportunity to him.  

The author is grateful to Michael Stoll for many useful discussions 
which contributed significantly to the improvement of the contents of this paper.


\begin{thebibliography}{99}

\bibitem[B]{Ber}
D.~J.~Bernstein, {\it Enumerating Solutions to p(a) + q(b) = r(c) + s(d)}, Mathematics of Computation,
Volume 70, Number 233, Pages 389-394

\bibitem[FMT]{FMT} J.~Franke, Y.~I.~Manin and Y.~Tschinkel, {\it Rational Points of Bounded Height on Fano Varieties},
Inventiones Mathematicae 95(1989), 421-435  

\bibitem[KaMa]{K&M} D.~Kanevsky, Yu.~I.~Manin. {\it  Composition of points and Mordell--Weil problem for
cubic surfaces.}
In: Rational Points on Algebraic Varieties
(ed. by E.~Peyre, Yu.~Tschinkel), Progress in Mathematics, vol. 199,
Birkh\"auser, Basel, 2001, 199--219. Preprint math.AG/0011198

\bibitem[Magma 1997]{Magma} {\sf MAGMA} is described in
  {\sc W. Bosma, J. Cannon} and {\sc C. Playoust:} 
  {\it The Magma algebra system I: The user language}, 
  J. Symb. Comp. {\bf 24}, 235--265 (1997).
  (Also see the Magma home page at 
  {\tt http://www.maths.usyd.edu.au:8000/u/magma/}\,.)

\bibitem[MA1]{Man1}
Yu.~I.~Manin. {\it Cubic Forms: Algebra, Geometry, Arithmetic.} North Holland, 1974 and 1986.

\bibitem[MA2]{Man}
Yu.~I.~Manin. {\it Mordell-Weil Problem for Cubic Surfaces.} In: Advances in the 
Mathematical Sciences---CRM's 25 Years (L. Vinet, ed.)  CRM Proc. and Lecture Notes, vol. 11, Amer. Math. Soc., 
Providence, RI, 1997, pp. 313--318. Preprint math.AG/9407009

\bibitem[PT]{P&T}
E.~Peyre and Y.~Tschinkel, {\it Tamagawa Numbers of Diagonal Cubic Surfaces of Higher Rank.
Rational points on algebraic varieties}, Progr. Math., 199, Birkh\"auser, Basel, 2001, pp. 275--305. 
Preprint math.AG/9809054

\bibitem[SlSD]{S&S-D} J.~B.~Slater and Sir P.~Swinnerton-Dyer, {\it 
Counting points in cubic surfaces, I}, Nombres et r\'epartitions
de points de hauteur born\'ee, Ast\'erisque 251 (1998), 1-11

\end{thebibliography}
\end{document}